\documentclass[11pt]{article}   	
\usepackage[top=.75in, bottom = .75in, left = .75in, right =.75in] {geometry}         		
\usepackage{graphicx}			
\usepackage{amssymb}
\usepackage{amsmath}
\usepackage{psfrag}
\usepackage{multicol}

\def\vp{\vspace{\baselineskip}}
\def\Nn{{\mathbb N}}
\parskip = \baselineskip
\begin{document}

\vp

\centerline{{\Large \sf A remark on Pascal's Triangle}}
\vp

\centerline{Chaim Goodman-Strauss}
\centerline{National Museum of Mathematics}
\centerline{chaimgoodmanstrauss@gmail.com}

\begin{abstract}
Through a series of elementary exercises, we explain the fractal structure of Pascal's triangle when written modulo $p$ using a weak form of an 1852 theorem due to  Kummer~\cite{ft}: A prime $p$ divides $\dfrac {n!}{i!j!} $ if and only if there is a carry in the addition $i+j=n$ when written in base $p$.
\end{abstract}

\vp
Centuries before Blaise Pascal (1623-1662), Pascal's Triangle was already known to  Yang Hui\raisebox{-.2\baselineskip}{\includegraphics[height=.9\baselineskip]{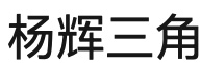}} \ (1238-1298), Jia Xian \raisebox{-.2\baselineskip}{\includegraphics[height=.9\baselineskip]{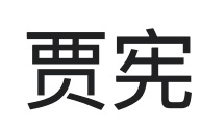}}\ (1010-1070) and 
Al-Karaji \raisebox{-.2\baselineskip}{\includegraphics[height=.85\baselineskip]{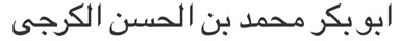}} \ (c. 953 – c. 1029).
These mathematicians and countless other since have shown how these  numbers are amazingly useful with many remarkable properties.\footnote{See Lecture 2 of D. Fuchs and S. Tabachnikov's wonderful {\em Mathematical Omnibus}~\cite{ft}, from which we learned of Kummer's theorem and its proof,  just one of many such riches in the lecture and book. In that lecture the authors give another proof of self-similarity of Pascal's triangle written modulo $p$.  See also~\cite{alex} for another take on  many of the same essential observations presented here. Please let me know if there are further sightings!}
  We are curious how a pattern like this one can appear just by marking which numbers in the triangle are even and which are odd:

\centerline{\includegraphics[width = 4.5in]{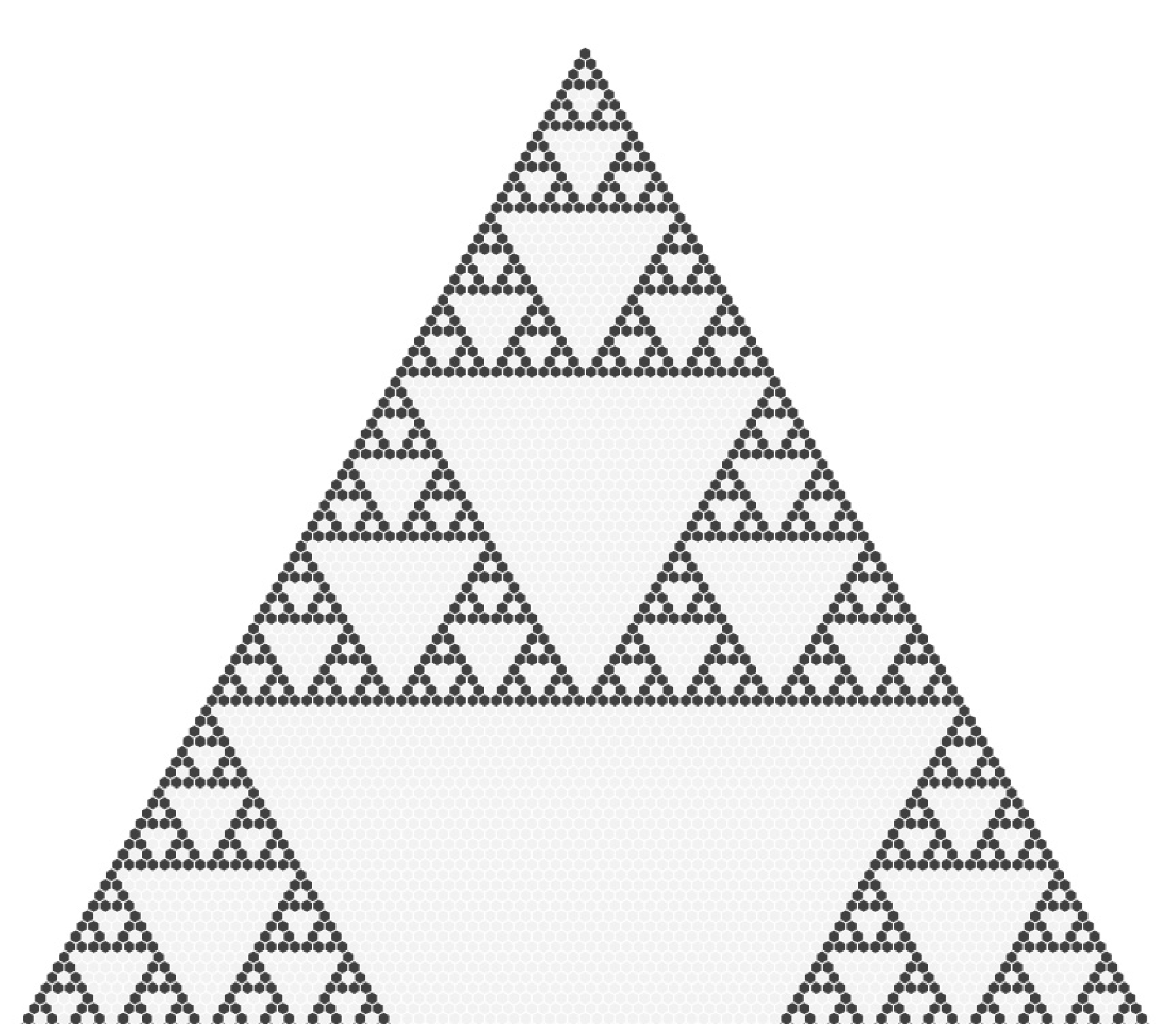}}
  
 \vp
You are likely familiar with the basic rule defining Pascal's Triangle: Beginning from a row of all 0's (which we ignore) and a single 1, the entries on each successive row are the sums of the neighbors above them.  We'll start another way around, defining coordinates on the triangle, and a formula for the value at each coordinate. The theorem below asserts that this formula satisfies the rule, and so really does give values for Pascal's Triangle --- you'll have to prove the theorem is true though!

We can measure the position of an entry by its coordinates $n$, $i$ and $j$, where  $n=i+j$ and   these are {\em natural numbers},\footnote{The {\em counting numbers} are $\{1,2,3,\dots\}$. The {\em integers}, usually denoted $\mathbb Z$, are defined as ${\mathbb Z}:=\{\dots, -3,-2,-1,0,1,2,3 \}$. The ``whole'' numbers usually are the counting numbers  additionally including 0, but we must always specify what we mean by  the ``natural'' numbers. Different authors and subcultures might  or might not consider 0 a natural number. Sometimes I do; sometimes I don't. I try to say which always.} that is, 
$n,i,j\in \Nn:=\{0,1,2,\dots\}$.
The row is measured by $n$; the  NW-SE diagonals are counted by $i$ and the NE-SW diagonals are counted by $j$.  For example, the shaded entry is on the $n=6$th row, on the $i=2$nd NW-SE diagonal, and the $j=4$th NE-SW diagonal.
If we define the entry at  $n,i$ and $j$ to be  $\dbinom n i := \dbinom n j := \dfrac{n!}{i!\ j!}$, these are the famous binomial coefficients, worthy of many long sets of notes.  Most importantly, $\dfrac{n!}{i!\ j!}$ is the coefficient of $x^iy^j$ in the expansion of $(x+y)^n$. This is also the number of ways to divide $n$ distinct things into two piles, an ``$x$'' pile with $i$ things in it, and a ``$y$'' pile with $j$ things in it, as well as the number of ways to choose $i$ out of $n$ things and ignore the other $j$ of them. 

{\bf Exercise} For  given $n,i,j\in\Nn$ with $n=i+j$, why are these numbers of ways to do things equivalent, and why does the formula correctly count them?

\vp {\bf Exercise} Look at the patterns in the coordinates along diagonals --- those running NW-to-SE and then those running NE-SW --- and along rows. 
Fill in the values of these entries in the first few rows. ($0!$ is defined to be equal to $1$. This makes sense for many reasons, not least of which that it makes sense for the values here.)

\scalebox{.7}{
\def\bbx#1#2#3{{\large $\dfrac{#1!}{#2!\ #3!}$}}
\psfrag{A}[cc]{\bbx n i j }
\psfrag{a}[cc]{\bbx000}
\psfrag{b}[cc]{\bbx110}
\psfrag{c}[cc]{\bbx101}
\psfrag{d}[cc]{\bbx220}
\psfrag{e}[cc]{\bbx211}
\psfrag{f}[cc]{\bbx202}
\psfrag{g}[cc]{\bbx330}
\psfrag{h}[cc]{\bbx321}
\psfrag{i}[cc]{\bbx312}
\psfrag{j}[cc]{\bbx303}
\psfrag{k}[cc]{\bbx440}
\psfrag{l}[cc]{\bbx431}
\psfrag{m}[cc]{\bbx422}
\psfrag{n}[cc]{\bbx413}
\psfrag{o}[cc]{\bbx404}
\psfrag{p}[cc]{\bbx550}
\psfrag{q}[cc]{\bbx541}
\psfrag{r}[cc]{\bbx532}
\psfrag{s}[cc]{\bbx523}
\psfrag{t}[cc]{\bbx514}
\psfrag{u}[cc]{\bbx505}
\psfrag{v}[cc]{\bbx660}
\psfrag{w}[cc]{\bbx651}
\psfrag{x}[cc]{\bbx642}
\psfrag{y}[cc]{\bbx633}
\psfrag{z}[cc]{\bbx624}
\psfrag{0}[cc]{\bbx615}
\psfrag{1}[cc]{\bbx606}
\psfrag{2}[cc]{\bbx770}
\psfrag{3}[cc]{\bbx761}
\psfrag{4}[cc]{\bbx752}
\psfrag{5}[cc]{\bbx743}
\psfrag{6}[cc]{\bbx734}
\psfrag{7}[cc]{\bbx725}
\psfrag{8}[cc]{\bbx716}
\psfrag{9}[cc]{\bbx707}

\centerline{\includegraphics[width=4.5in]{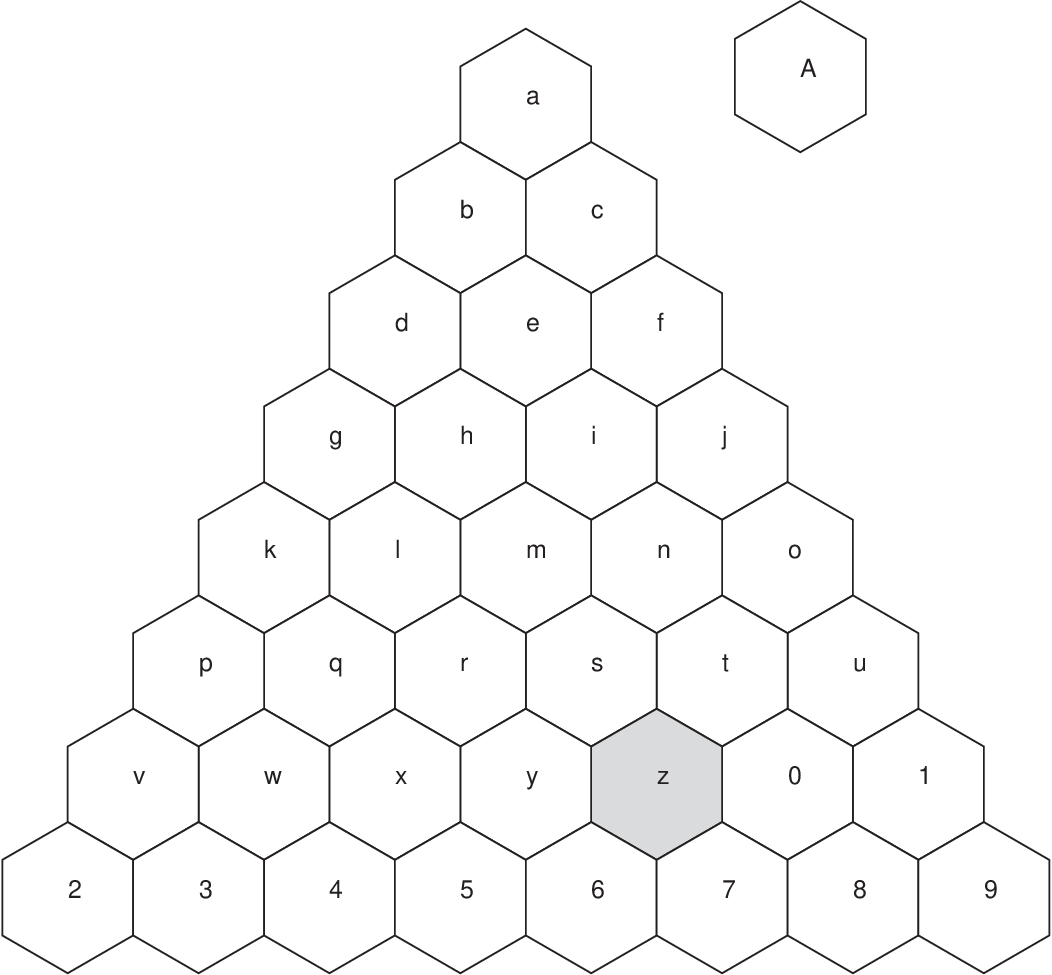}}}

Check that a cell at position $n,i,j$ has upstairs neighbors at $(n-1),(i-1),j$ and $(n-1),(i-1),j$, and prove that the binomial coefficients satisfy the rule defining Pascal's Triangle

\vp
{\bf Theorem 1} {\em For all $n,i,j\in \Nn$ with $n=i+j>1$, $\dbinom {n} {i} = \dbinom {n-1} {i} + \dbinom {n-1} {i-1}$.
}

In other words, show  $$\frac{n!}{i!\ j!}  = \frac{(n-1)!}{i!\ (j-1)!}+\frac{(n-1)!}{(i-1)!\ j!}$$

To see what is going on, test a few examples, such as $$\frac{8!}{5!3!}=\frac{7!}{5!2!}+\frac{7!}{4!3!}\ \text{ \ or }$$
$$\frac{32!}{19!13!}=\frac{31!}{19!12!}+\frac{31!}{18!13!}.$$

\vp\vp
 This next theorem may seem obvious, but why is it true -- why should {\em all} of the factors of $i!$ and $j!$ cancel out with factors of $n!$\ ?

{\bf Theorem 2} {\em For all $n,i,j\in\Nn$ with $n=i+j$,  $\dbinom {n} {i}$ is an integer.}

If you've completed the exercise above, then you provide a simple, but indirect proof of  Theorem 2: $n!/i!j!$ counts  whole numbers of things. But in a more hands on way, how might you   begin from the assumption that $\binom{0}{0}=1$ and $\binom 0 k =0$ otherwise, and then use Theorem 1.

\vp\vp\vp\vp

We are interested in Pascal's Triangle modulo $p$, for  primes $p$, looking at the remainders of the entries when divided by $p$. For example, we could color the entries   $\dbinom n {i}$  ``white'' if they are even (are 0 mod 2) or ``black'' if they are odd (are 1 mod 2). Equivalently, we could generate Pascal's Triangle, finding each entry by adding the values of the entries above, but only working modulo $p$. 

Modulo 2, we only care whether numbers are {\Large $e$}ven or od{\Large  $d$}. We have that $e+e=d+d=e$ and $e+d=d+e=d$. 

This is a kind of ``cellular automaton''. Each cell depends on the states of the  cells above it, and we can think of this as a movie of a row of cells changing over time:  Each row is the movie frame after the row above it.

\vp\vp

{\bf Exercise} Fill in a couple of dozen rows of Pascal's Triangle mod 2 to see a pattern! (Or turn to figures on the next few pages.) 

{\bf Exercise} Fill in a couple of dozen rows of Pascal's Triangle mod 3, mod 5 etc, or write a program to do so. What is the rule for how these patterns are generated?

\vp

\centerline{\includegraphics{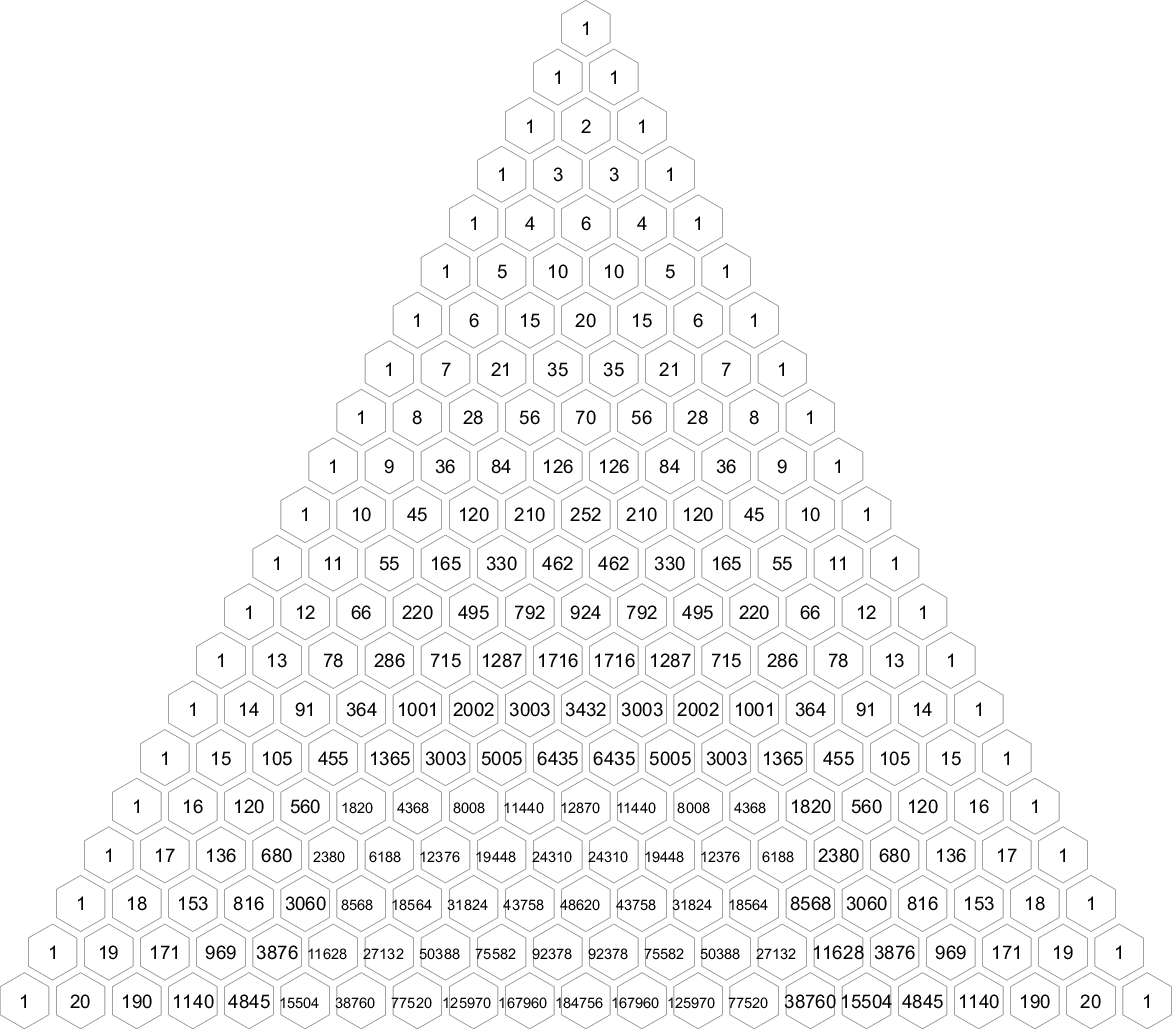}}

\newpage

\newpage

\centerline{
\includegraphics[width=\textwidth]{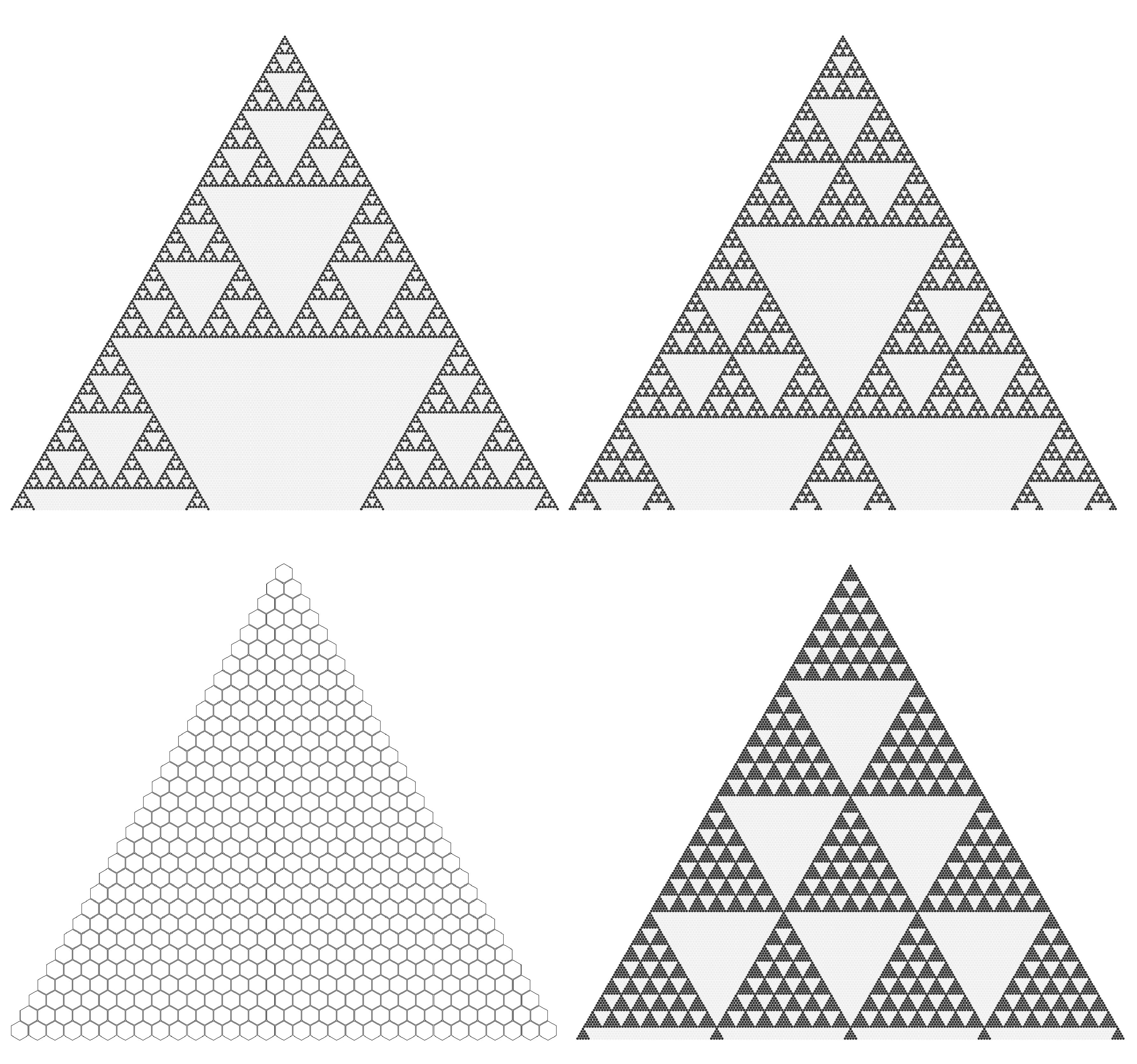}}

Two-hundred rows of cells are colored white if they are divisible by 2 (top left) or 3 (top right) or 7 (bottom right), and colored black otherwise.  
At bottom left, you can shade in thirty rows of cells mod 5: What should the pattern be? 

(In the bottom right figure how many rows are needed before the next-sized white triangle appears?)

(What will the pattern be for composite numbers? For powers of primes?)

\newpage

{\sf\Large Place-value representations:} 

\vp
We are familiar with decimal representations such as $5342= 5\cdot 10^3 + 3\cdot 10^2+4\cdot 10^1+2\cdot 10^0$, and  for any integer $b>1$ we can express counting numbers base $b$, as a sum of multiples of powers of $b$. For example, in base 7, $5342_7= 5\cdot 7^3 + 3\cdot 7^2+4\cdot 7^1+2\cdot 7^0 = 1892_{10}$. Of course for base $b$, we only need  $0,1,...,(b-1)$ as digits.

There is a simple procedure to find the digits of a number base $b$: Repeatedly divide by $b$; the remainders are the digits written in reverse order. For example:

$\begin{array}{rl} 2932& = 325\cdot 9+\textcircled{7}\\
325 & = 36\cdot 9 + \textcircled{1}\\
36&=4\cdot 9 + \textcircled{0}\\
4&=0\cdot 9 + \textcircled{4}\\
\end{array}$

Reading from bottom to top, this shows us that  $2932_{10}=4017_9$ (which you can verify).
This works because $2932 = (((0\cdot 9 + 4)\cdot 9+0)\cdot 9 + 1)\cdot 9+7 = 4\cdot 9^3 + 0\cdot 9^2 + 1\cdot 9 + 7$. 

{\bf Exercise} Try out a zillion examples, in all kinds of bases. Check your work, converting back into base 10. 

For a number $n$ written in base $b$, we'll call the digits of $n$ as ... $n_2$, $n_1$, $n_0$, which take values from $0$ to $b-1$. That is, $n = ... + n_2 b^2+n_1 b^1 + n_0$.

{\bf A lemma} {\em In base $b\in\Nn$, the $k$th digit $n_k$ of $n$  is $$n_k=\left\lfloor n/b^k\right\rfloor - b\cdot \left\lfloor n/b^{k+1}\right\rfloor$$} where $\left\lfloor x \right\rfloor$ is the floor function, returning the greatest integer not larger than $x$. 

Try out some examples to see if you agree with this, and get some intuition. Then supply a proof. \vp

Here is a nice theorem. Experiment with examples to see how it works. We'll prove it in a few steps.

{\bf Theorem 3} {\em Let $p$ be a prime number and let $n$ be an integer. Let $N$ be the sum of the digits of $n$ written in base $p$. Then the number of factors of $p$ in $n!$ is exactly $$(n-N)/(p-1).$$}

{\bf Exercise} 
As a check how many factors of 5 are there in 132!\ ? What is 132 written in base 5? For both of these questions, you would like to know, from 1 to 132,  how many multiples of 5 there are, how many multiples of $5^2$ there are, how many multiples of $5^3$ there are, etc.

{\bf Exercise} How many factors of 7 are there in 365! ? Make up and work out more examples!

{\bf Another lemma}  {\em The number of factors of a prime $p$ in $n!$ is $$\left\lfloor  n/ p \right\rfloor+ 
\left\lfloor  n /{p^2} \right\rfloor+ \left\lfloor  n/ {p^3} \right\rfloor+\dots .$$ }
 The terms of this series are eventually 0. As always, try several examples. (From~\cite{alex} we learned that  Adrien-Marie Legendre presented this  lemma in 1808.)
\vp

Why is the lemma true? Suppose that written in base $p$, $n= b_k b_{k-1} \dots b_1 b_0$. (In other words, there are $k+1$ digits,  $b_k, b_{k-1}, \dots, b_1, b_0$ and $n=b_k p^k + \dots b_1 p^1 + b_0$.) 

How many numbers from 1 to $n$ are divisible by $p$? (each contributing a factor of $p$ to $n!$)

How many numbers from 1 to $n$ are divisible by $p^2$? (each contributing another factor of $p$ to $n!$.) Etc. 

The largest power of $p$ that divides $n!$ is $p^k$. How many numbers from 1 to $n$ are divisible by $p^k$? (each contributing another factor of $p$ to $n!$.)

All together, what is the relationship between the digits of $n$ base $p$ and the number of factors of $p$ in $n!$ ? (Did we use that $p$ is prime?)

\vp

Let's put our two lemmas together and prove  Theorem 3. Let $N$ be the sum of the digits of $n$ written in base $p$, a prime. The first lemma gives us a formula for each digit, which we can use to produce a formula for $N$.

$$\begin{array}{rl} N &= \sum_{k=0}^{\infty} \left\lfloor \frac n {p^k}\right\rfloor - p\cdot \left\lfloor \frac n {p^{k+1}}\right\rfloor\\\\
& = \left\lfloor \frac n {p^0}\right\rfloor+ (1-p) \sum_{k=1}^{\infty} \left\lfloor \frac n {p^k}\right\rfloor\\\\
& = n +  (1-p) f
\end{array} $$ where $f$ is the number of factors of $p$ in $n!$. Solving for $f$, we have the theorem.

\vp\vp

{\sf\Large Addition}

Remind yourself how to add, especially in another base. In base 7, for example, our digits are 0,1,2,3,4,5 and 6, and we carry whenever our total is bigger than 6. Add this base 7:

\vp\vp
$\begin{array}{rr}
&253_7\\
+&415_7\\\hline
\\
\end{array}$

\vp
Hopefully you have $1001_7$ as your answer, with three carries. (You can, and should, check your work by converting your work to base 10.) 

\vp

{\bf Exercise} Try out some more addition problems in different bases, and check your work base 10. 

\vp How does addition change the {\em total } of the  digits? In the example $253_7+415_7=1001_7$, we begin with a total of $2+5+3+4+1+5$ in our digits, but end up with only $1+1$, a decrease of 18. Why 18? 

Each carry takes a total of 7 from the digits and replaces it with a 1, decreasing the total by 6, and there are three carries. 

\vp
More generally, working base $b$,  how many carries when we  add $i+j$ to get $n$? Can you prove this theorem?

\vp
{\bf Theorem 4} {\em  Let $N$, $I$ and $J$ be the sums of the digits of $n$, $i$ and $j$ written in base $b$. Let $c$ be the number of carries of $n=i+j$ in base $b$. Then $c = (I+J-N)/(b-1)$.
}

\vp

Is it obvious, or not, that this formula for $c$ is never negative? (Why is $(I+J-N)$ always non-negative? Why is the sum of digits of $I$ and $J$ always at least  that of $N$?)

Now you can put these pieces together to prove:

\vp

{\bf Theorem} (Kummer, 1852) {\em Let $p$ be prime and $i,j$ and $n=i+j$ be natural numbers. The number of factors of $p$ in $\dfrac{n!}{i!\; j!}$ is the number of carries in $i+j=n$ written in base $p$.}

\vp

Where and how do we use that $p$ is prime in this theorem? Can you find a  correct theorem for composite numbers $p$? 
\vp

{\sf\Large Carries in Pascal's Triangle}

How do we get the fractal pattern in Pascal's Triangle modulo 2, or modulo any prime $p$? 

Let's work in binary from this point forward, but keep your eye on how this will work for any prime base. 

We want to understand which cells are to be colored white--- which values of $n!/i!j!$ are even. Kummer's theorem tells us these are precisely those that have at least one carry in the binary addition of $i+j=n$. 

Which binary additions have a carry? The sum $1+1$  equals $0$, carry 1. The only other carries are if we are already carrying something else. For example, add this:

\vp 
$\Large \begin{array}{rr}
&1011\\
+&1001\\\hline\\\end{array}$

The right most sum 1+1 yields 0, carrying a 1. In the 2's place, we have 1+0, plus the carried 1, yielding 0, again carry 1. This carrying stops, because we now are adding 0+0, plus the carried 1, yielding 1 with no carry.  
Let's call a carry a {\em stopping carry} if it stops, if there is not a carry to in the place-position to the left. So the carry in the 4's position is a stopping carry in the example we just saw. 

In fact

{\bf Lemma} {\em Suppose there is a carry in a binary sum. Then  in the place position of any  stopping carry, the sum will be $0+0 = 1$ with a carried $1$.}

To prove this lemma, a stepping stone kind of theorem, just check: If there is a carry at all, could any {\em other} sum be a stopping carry? 

Since every carry stops, every sum that has a carry has a place position with values $0$ for $i$ and $j$ and $1$ for $n$. Let's call these place positions ``special'' and these will be the positions we will pay attention to, because if {\em any} place position is ``special'' for $n$, $i$ and $j$, then the cell in that position in Pascal's Triangle will be colored white. 

{\bf Exercise} Which sums have carries, and where are the stopping carries? Which place positions have values $0$ for $i$ and $j$ and $1$ for $n$ and are ``special''. 

\hspace*{-2em}
\begin{minipage}{\textwidth}
\begin{multicols}4
\begin{itemize}
\item $\begin{array}{rr}&10010\\+&1100 \\\hline\\ \end{array}$
\item $\begin{array}{rr}&111\\+&1\\\hline\\ \end{array}$
\item $\begin{array}{rr}&1101\\+&1001\\\hline\\ \end{array}$
\item $\begin{array}{rr}&10101010\\+&10010010\\\hline\\ \end{array}$
\end{itemize}
\end{multicols}
\end{minipage}

\vp

{\bf Exercise} In binary, in Pascal's Triangle, label 
the rows and the two families of diagonals. 

\vp
\centerline{
\includegraphics[scale = .7]{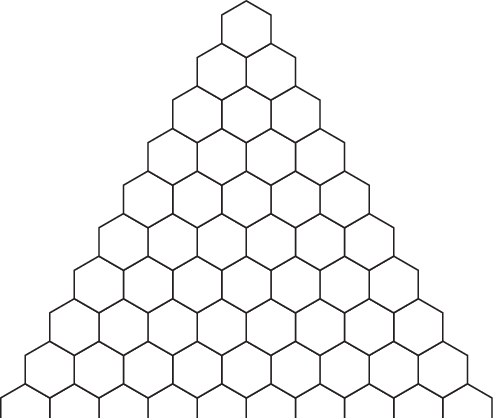}}

We are interested in the pattern of carries in Pascal's Triangle mod 2.   If there are any carries, there are special places  and so we will consider those. 

\vp

{\bf Exercise} What cells could have their 1's place be special? 

This is a trick question, because it is impossible for $i$ and $j$ to be $0$ in the 1's place and  $n$ to be $1$ in the 1's place.

\vp
{\bf Exercise} What cells in could be have their 2's place be special? What fraction of cells are special in this place?

\vp

Let's color this in, shading in the cells that {\em aren't} special. First, $n$ must have a $1$ in its 2's place: In binary, $n$ must be of the form $\dots**1*$. Shading out the cells that aren't of this form will produce some horizontal stripes, which you can shade in here:

\includegraphics{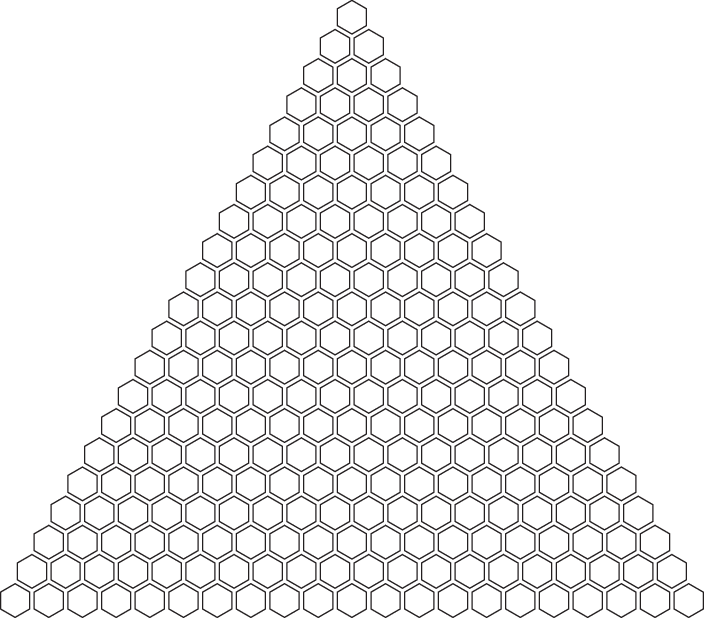}

Similarly, you can shade out the cells that aren't of the right form for $i$ and $j$. These must have a $0$ in their 2's place for the cell to be special in the 2's place, and if {\em either} of them don't the cell {\em isn't}. Shading out the cells that don't have a $0$ in the 2's place of $i$ produces diagonal stripes, marching NW-to-SE. Shading out the cells that don't have a $0$ in the 2's place of $j$ produces diagonal stripes, marching NE-to-SW.  Only the cells that are not in any of these shaded out stripes, will satisfy all three conditions on the 2's place. 

\vp

\centerline{
\includegraphics[width = 2.5in]{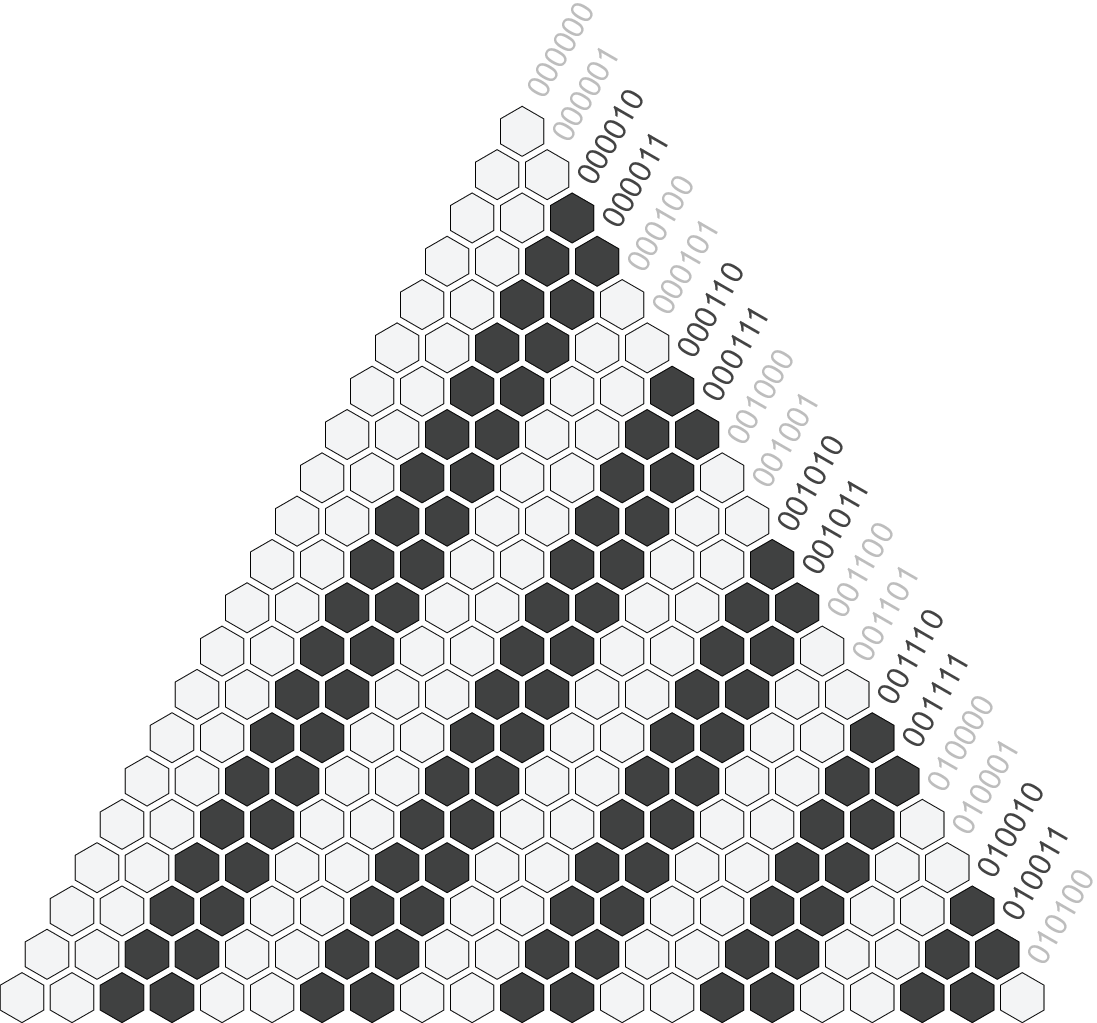}
\includegraphics[width = 2.5in]{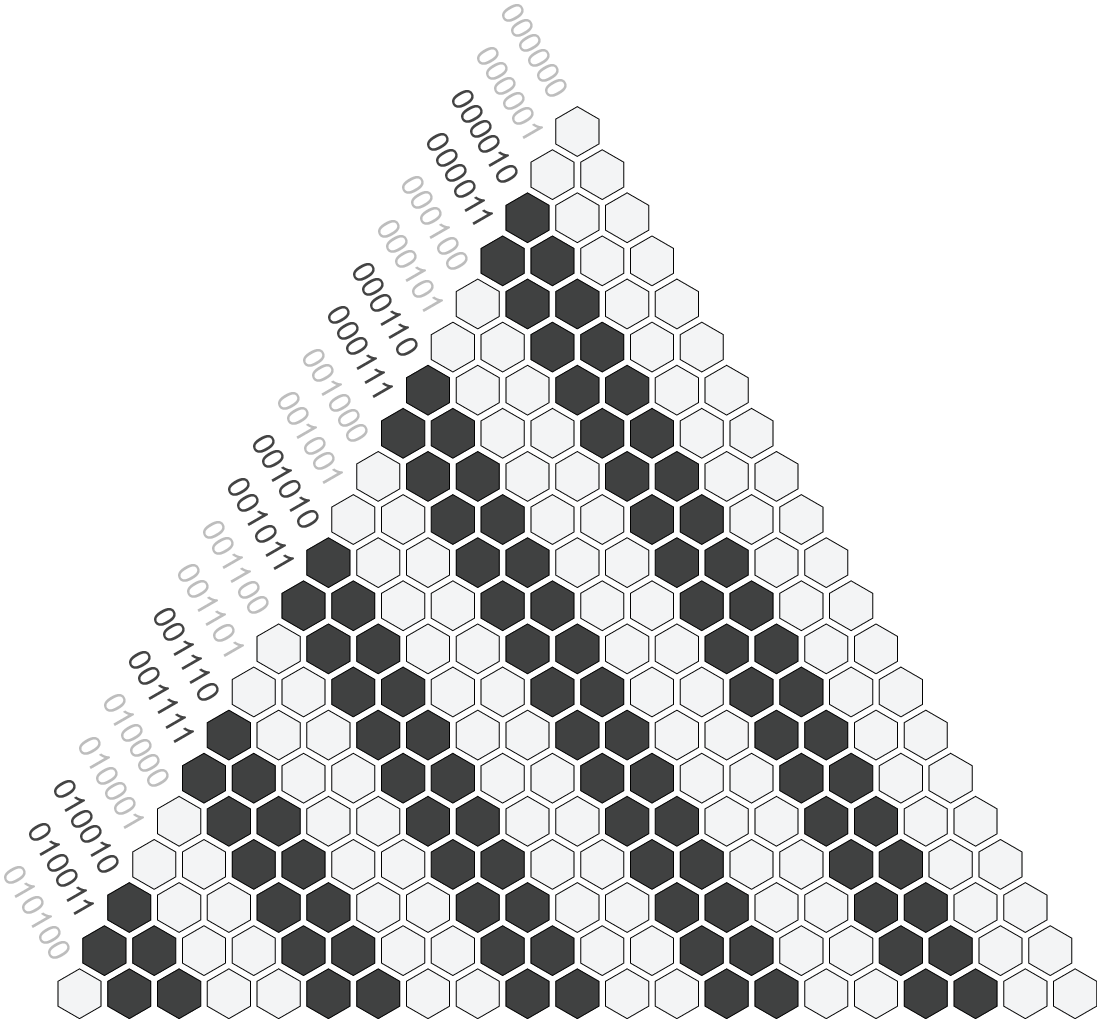}
\includegraphics[width = 2.5in]{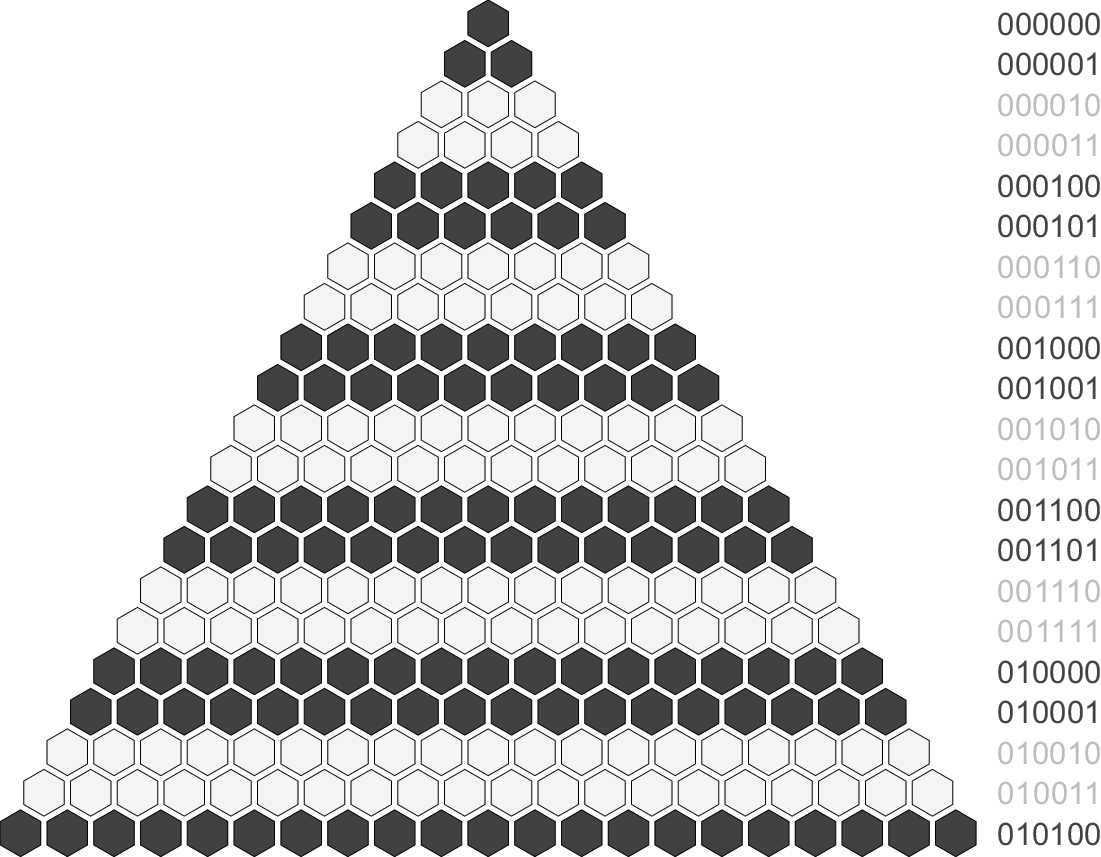}}

\vp
{\bf Exercise} Which cells in Pascal's Triangle are special in the  4's place?  (Where and how wide are the stripes that eliminate special cells?) What fraction of cells are special in the 4's place?

{\bf Exercise} Which cells in Pascal's Triangle are special in the 8's place? The 16's? Etc?

\includegraphics[width = \textwidth]{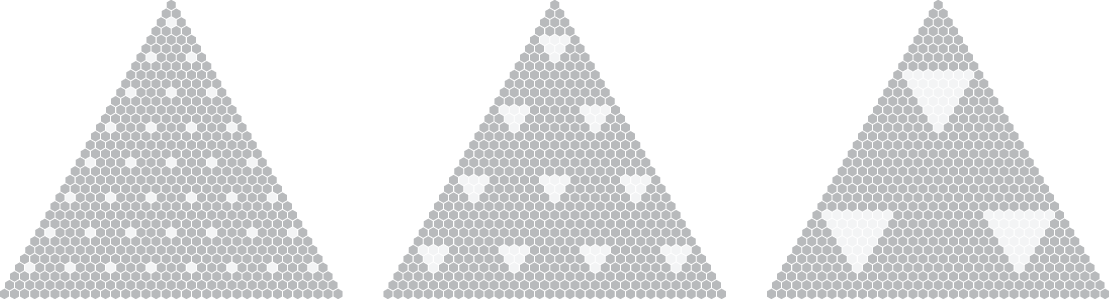}

\vp
{\bf Exercise} Which cells in Pascal's Triangle are special in {\em some} place?

\vp We can now describe the fractal structure in Pascal's triangle modulo 2:  It is the union of  lattices of triangles, each lattice being the intersection of stripes. These correspond to the locations $n!/i!j!$ in the Pascals triangle in which, for some place position, $i$ and $j$ have a bit $0$ and $n$ has a bit $1$, in other words for the locations in which there is a carry in the binary addition $i+j=n$. In turn, these are precisely the values $n!/i!j!$ that are even.

{\bf Exercise} How does this generalize for other primes?

\end{document}